\documentclass[letterpaper,11pt]{amsart}

\usepackage[margin=1.2in]{geometry}
\usepackage{amsmath,amsthm,amssymb, mathtools}
\usepackage{xspace,xcolor}
\usepackage[breaklinks,colorlinks,citecolor=teal,linkcolor=teal,urlcolor=teal,pagebackref,hyperindex]{hyperref}
\usepackage[alphabetic]{amsrefs}
\usepackage{comment}
\usepackage{enumerate}
\usepackage[all]{xy}
\usepackage{tikz-cd}
\usepackage{color}
\usepackage{comment}

\theoremstyle{plain}
\newtheorem{theo}{Theorem}[section]

\newtheorem{lemm}[theo]{Lemma}
\newtheorem{exam}[theo]{Example}

\newtheorem{coro}[theo]{Corollary}
\newtheorem{conj}[theo]{Conjecture}
\newtheorem{ques}[theo]{Question}

\theoremstyle{definition}

\theoremstyle{remark}
\newtheorem{rema}[theo]{Remark}

\newcommand{\eps}{\epsilon}

\newcommand{\PP}{\mathbb{P}}
\newcommand{\QQ}{\mathbf{Q}}

\newcommand{\RR}{\mathbf{R}}

\newcommand{\VV}{\mathbb{V}}

\newcommand{\cF}{\mathcal{F}}

\newcommand{\cO}{\mathcal{O}}

\newcommand{\uind}[1]{^{(#1)}}

\newcommand{\lsta}{_{*}}

\newcommand{\sta}{^{*}}

\newcommand{\lround}[1]{\left\lfloor #1 \right\rfloor}

\DeclareMathOperator{\rank}{rank}

\DeclareMathOperator{\codim}{codim}

\DeclareMathOperator{\SheafHom}{\mathcal{H}om}

\DeclareMathOperator{\ord}{ord}

\DeclareMathOperator{\supp}{supp}

\usepackage{soul}

\begin{document}
	\thanks{The author was partially supported by NSF grant DMS-2301463 and the Simons Collaboration grant Moduli of Varieties}
	
	\subjclass[2020]{14D06, 14E30, 14J15, 32J25}
	
	\author{Hyunsuk Kim}
	
	\address{Department of Mathematics, University of Michigan, 530 Church Street, Ann Arbor, MI 48109, USA}
	
	\email{kimhysuk@umich.edu}	

	\begin{abstract}
		We observe what the canonical bundle formula gives towards a conjecture of Schnell on algebraic fiber spaces, a question concerning the equivalence between the non-vanishing conjecture and the Campana--Peternell conjecture. As a result, we give a partial result on Schnell's conjecture under two independent assumptions. One weakens Schnell's assumption of the pseudo-effectivity of the canonical bundle of the base by adding some effective divisor supported on the ramification locus. The other is analogous to results on algebraic fiber spaces where the existence of good minimal models of a general fiber is assumed, but we use a priori a weaker assumption. More precisely, we prove Schnell's conjecture when the canonical class of the general fiber is represented by a \textit{rigid current}.
	\end{abstract}
	
	\title[Canonical bundle formula and a conjecture by Schnell]{Canonical bundle formula and a conjecture on certain algebraic fiber spaces by Schnell}
	
	\maketitle
	
	
	
	\section{Introduction}
	
	In \cite{Schnell-conj}, Schnell attempts to prove an equivalence between two conjectures, namely the non-vanishing conjecture and the Campana--Peternell conjecture (on Kodaira dimension, see \cite{CP-conj}*{page 43}). Both conjectures are special cases of the abundance conjecture. The non-vanishing conjecture predicts that if $X$ is a smooth projective variety whose canonical bundle $K_{X}$ is pseudo-effective, then a positive multiple of $K_{X}$ is effective. The Campana--Peternell conjecture says the following.
	
	\begin{conj}[Campana--Peternell] \label{conj:Campana-Peternell}
		Let $X$ be a smooth projective variety and let $D$ be an effective divisor. Suppose that $m_{0}K_{X} - D$ is pseudo-effective for some $m_{0} > 0$. Then $\kappa(X) \geq \kappa(X, D)$.
	\end{conj}
	In particular, if we put $D = 0$ in the conjecture above, we get the non-vanishing conjecture. In this sense, the Campana--Peternell conjecture seems to be a priori stronger than non-vanishing. In \cite{Schnell-conj}, it is shown that assuming non-vanishing, the Campana--Peternell conjecture reduces to proving the following problem on algebraic fiber spaces.
	
	\begin{conj} \label{conj:alg-fibr-space}
		Let $f \colon X\to Y$ be an algebraic fiber space between smooth projective varieties with $\kappa(F) = 0$, where $F$ is the general fiber. If there is an ample divisor $H$ on $Y$ such that $m_{0}K_{X} - f\sta H$ is pseudo-effective, then $\kappa(X) = \dim Y$.
	\end{conj}
	
	This formulation is slightly different than \cite{Schnell-conj}*{Conjecture 10.1}. However, one can reduce to the case above by following the lines of \cite{Schnell-conj}*{\S9}. We prefer this formulation so that the canonical bundle formula naturally comes into the picture. Schnell proves his conjecture when $K_{Y}$ is pseudo-effective, using the result of P\u{a}un--Takayama \cite{Paun-Takayama} on semi-positive metrics on relative pluri-adjoint line bundles.
	
	The aim of this article is to observe what canonical bundle formula gives toward Conjecture \ref{conj:alg-fibr-space}. As a result, we are able to weaken Schnell's assumption on the pseudo-effectivity of $K_{Y}$, and also observe that \textit{rigid currents} in the sense of \cite{lazic2024rigid} naturally plays a role in this conjecture.
	
	\subsection{Canonical bundle formula} \label{subsec:CBF}
	We introduce the canonical bundle formula, following \cite{Kaw-subadjunctionII, FM-canonical-bundle-formula, Kollar-KodCBF}, in the set-up of Conjecture \ref{conj:alg-fibr-space}. We also refer to \cite{danokim-2019canonical} for analytic aspects of the canonical bundle formula using fiberwise integrals, which is related to the content appearing in the later part of this paper.
	
	First, following \cite{Schnell-conj}*{\S5}, we note that this conjecture is insensitive upon taking birational models of $X$ and $Y$. By taking suitable resolutions, we may and will assume that we have an SNC divisor $Q = \sum_{l \in L} Q_{l}$ on $Y$ such that $f$ is smooth over $Y_{0} = Y \setminus Q$, and $f^{-1}(Q)_{\mathrm{red}} = \sum_{j \in J} P_{j}$ is an SNC divisor on $X$. We can write $K_{X}$ in the following way:
	\begin{equation} \label{equation:CBF}
		K_{X} \sim_{\QQ} f\sta (K_{Y} + B_{Y} + M_{Y}) + \Delta, 
	\end{equation}
	so that
	\begin{enumerate}
		\item $f\lsta \cO_{X}(\lround{i\Delta^{+}}) = \cO_{Y}$ for all $i > 0$,
		\item $\codim_{Y} f(\supp \Delta^{-}) \geq 2$, and
		\item $H^{0}(X, mK_{X}) = H^{0}(Y, m(K_{Y} + B_{Y} + M_{Y}))$ for sufficiently divisible $m$ (see \cite{FM-canonical-bundle-formula}*{Theorem 4.5}).
	\end{enumerate}
	
	Here, $\Delta^{+}$ and $\Delta^{-}$ are the effective divisors sharing no irreducible components and satisfying $\Delta = \Delta^{+} - \Delta^{-}$. Roughly speaking, the divisor $B_{Y}$ measures the singularity of the fibers of $f$, and $M_{Y}$ measures the Hodge-theoretic variation of the fibers. In particular, $B_{Y}$ is always effective and is supported in the discriminant locus of $f$, and $M_{Y}$ is nef. For a precise description of the divisors $B_{Y}$ and $\Delta$, we refer to \S\ref{subsec:coeff-CBF}.
	
	We return to the setup of Conjecture \ref{conj:alg-fibr-space}. Suppose we have some $m_{0} >0$ such that $m_{0}K_{X} - f\sta H$ is pseudo-effective. Hence $K_{X} - \eps f\sta H$ is pseudo-effective for $\eps > 0$ small enough. By taking birational models of both $X$ and $Y$, we assume that $f \colon X \to Y$ satisfies the SNC assumption as above. In order to show that $\kappa(X) = \dim Y$, it is equivalent to showing that $K_{Y} + B_{Y} + M_{Y}$ is big, as the pluricanonical sections of $X$ descend to $Y$. Hence, it is enough to show that $K_{Y} + B_{Y} + M_{Y} - \eps H$ is pseudo-effective for some $\eps > 0$. Rewriting the canonical bundle formula, we have
	$$ (K_{X} + \Delta^{-} - \eps f\sta H) - \Delta^{+} \sim_{\QQ} f\sta (K_Y + B_Y +M_Y - \eps H).$$
	Note that $K_X + \Delta^- - \eps f\sta H$ has a semi-positively curved singular hermitian metric on $X$. Hence, we see that $\Delta^{+}$ is the only obstacle for descending the metric positivity of $K_{X} + \Delta^{-} - \eps f\sta H$ down to $Y$. In our main result, we give two separate ways to deal with this $\Delta^{+}$.
	
	\subsection{Main results}
	We state the main results.
	\begin{theo} \label{theo:main-theorem}
		As in the setup of Conjecture \ref{conj:alg-fibr-space} and \S\ref{subsec:CBF}, suppose that one of the two following conditions holds:
		\begin{enumerate}
			\item $K_{Y} + (1-\eps)B_{Y}$ is pseudo-effective for some $\eps > 0$, or
			\item for the general fiber $F$, the class $\{K_{F}\}$ is rigid.
		\end{enumerate}
		Then $\kappa(X) = \dim Y$.
	\end{theo}
	
	For a compact Kähler manifold $X$, a class $\alpha \in H^{1,1}(X) (\simeq H_{\mathrm{BC}}^{1,1}(X))$ is \textit{rigid} if it contains exactly one closed positive (1,1)-current $T$ such that $\{T\} \in \alpha$. For varieties of Kodaira dimension zero, the rigidity of the canonical class is tightly related to the existence of good minimal models. For a systematic treatment of rigid currents in terms of birational geometry, we refer to \cite{lazic2024rigid}. Here, we mention the easy direction that if $F$ admits a good minimal model, then the cohomology class $\{K_{F}\}$ is rigid \cite{lazic2024rigid}*{Lemma 3.9}.
	
	We point out that the proofs for the two cases are independent. Also, due to technical reasons concerning the volume asymptotics, we were not able to handle the case when $\eps = 0$ for the first condition.
	
	\subsection{Comparison with known results}
	We compare Theorem \ref{theo:main-theorem} with previously known results. The first result can be understood as a strengthening of \cite{Schnell-conj}*{Theorem 12.1}, where he assumes that $K_{Y}$ is pseudo-effective. Indeed, adding an effective divisor $(1-\eps)B_{Y}$ in the first case of Theorem \ref{theo:main-theorem} can be thought of as a correction term coming from the singularity of the fibers of $f$. Our initial motivation was to overcome the non-pseudo-effectivity of $K_{Y}$ by considering the contribution from the singular fibers, since heuristically, a projective morphism from a variety with `positive' canonical bundle to a variety with `negative' canonical bundle tends to have a lot of singular fibers (see \cite{Catanese-Schneider:polybounds}*{Question 4.1}, \cite{Kovacs}, \cite{Viehweg-Zuo}, \cite{Campana-Paun},\cite{park2022logarithmic} for related results). However, $K_{Y} + (1-\eps) B_{Y}$ is not always pseudo-effective, especially when the singular fibers are mildly singular, as explained in the following example pointed out by Niklas Müller.
	
	\begin{exam}
		Take a general elliptic K3 surface $f : X\to \PP^{1}$. Then $f$ has 24 singular fibers with at most nodal singularities. Consider $\nu :C \to \PP^{1}$ a double cover ramified over two points avoiding the singular values (we have $C \simeq \PP^{1}$). Then we obtain $g: Y \to C$ by taking the base change of $f$ by $\nu$, as
		$$ \begin{tikzcd}
			Y \ar[r,"\mu"] \ar[d, "g"] & X \ar[d, "f"] \\ C \ar[r, "\nu"] & \PP^{1}.
		\end{tikzcd}$$
		Using the ramification formula, we have
		$$K_{Y} \sim g\sta \cO_{C}(2).$$
		Hence, $K_{Y} - g\sta \cO_{C}(1)$ is pseudo-effective. However, if we write the canonical bundle formula
		$$ K_{Y} \sim_{\QQ} g\sta (K_{C} + B_{C} + M_{C}),$$
		we have $B_{C} = 0$ and $M_{C} \simeq \cO_{C}(4)$. Hence $K_{C} + (1-\eps)B_{C} = K_{C}$ is non-pseudo-effective for any $\eps > 0$.
	\end{exam}

	Our second result is somewhat analogous to results on algebraic fiber spaces assuming that the general fibers have good minimal models. However, we emphasize that the rigidity of the canonical class is a property \textit{purely in terms of hermitian metrics}. For example, a numerically trivial line bundle is rigid, even if it is non-torsion. This condition is weaker than the existence of good minimal models as mentioned in the section above, but still stronger than the Campana--Peternell conjecture. Indeed, given a smooth projective variety $X$ whose canonical class is rigid (which implies $\kappa(X) = 0$), if one writes $K_{X} = D + E$ where $D$ is an effective $\QQ$-divisor and $E$ is pseudo-effective, the rigidity of the canonical class forces $\kappa(X, D) = 0$. In this vein, we would like to formulate the following two questions:
	
	\begin{ques}
		Is there an example of a smooth projective variety $X$ with zero Kodaira dimension, where we know whose canonical class $\{ K_{X}\}$ is rigid, but the existence of good minimal models is not known?
	\end{ques}

    \begin{ques}
        Let $X$ be a smooth projective variety with zero Kodaira dimension. Suppose that if we write $K_{X} = D + E$ as a sum of an effective divisor $D$ and a pseudo-effective divisor $E$, we have $\kappa(X, D) = 0$. In this case, is $K_{X}$ rigid?
    \end{ques}
	
	We also would like to mention that assuming that the general fibers admit a good minimal model, Conjecture \ref{conj:alg-fibr-space} is known to be proved using purely algebraic methods. Here, we sketch the proof. Since the general fiber has a good minimal model, $X$ has a good minimal model $g \colon X' \to Y$ over $Y$. This is a special case of \cite{Hacon-Xu-lc-closures}*{Theorem 2.12}. Then by \cite{Taji-birational-positivity-dim-4}*{Lemma 3.4}, we have a birational modification of $g$ as the following diagram:
	$$ \begin{tikzcd}
		X \ar[r, dashed] \ar[rd, "f"'] & X' \ar[d, "g"] & \widetilde{X} \ar[l, "\mu"'] \ar[d, "\widetilde{g}"] \\ & Y & \widetilde{Y} \ar[l, "\pi"'],
	\end{tikzcd}$$
	so that $\mu\sta K_{X'} = \widetilde{g}\sta G$ for some $\QQ$-Cartier divisor $G$ on $\widetilde{Y}$. Note that $m_{0}K_{X'} - g\sta H$ is pseudo-effective, hence $\widetilde{g}\sta (m_{0}G - \pi\sta H)$ is pseudo-effective. Therefore, $m_{0}G - \pi\sta H$ is pseudo-effective. Since $\pi\sta H$ is big, $G$ is also big, and we see that $\kappa(X) = \kappa(X') = \kappa(\widetilde{Y}, G) = \dim Y$.
	
	We would also like to compare our second result with the inductive approach by \cite{CP-conj}*{Proposition 2.6}. Let $X$ be a smooth projective $n$-fold and let $D$ be an effective divisor such that $m_{0}K_{X} - D$ is pseudo-effective for some $m_{0}>0$. Assuming the Campana--Peternell conjecture in dimension $\leq n-1$, by applying the easy addition (see \cite{Mori-Bowdoin}*{Corollary 1.7}) to the Iitaka fibration of $X$, they show that the Campana--Peternell conjecture holds in dimension $n$ except for the two following situations:
		\begin{enumerate}
			\item $\kappa(X) = 0$ and $\kappa(X, D) >0$, or
			\item $\kappa(X)= -\infty$ and $\kappa(X, D) \geq 0$.
		\end{enumerate}
	Assuming the non-vanishing conjecture in dimension $n$ automatically rules out the second possibility. The point of the second statement of Theorem \ref{theo:main-theorem} is that assuming the rigidity of the canonical class of varieties of Kodaira dimension zero for dimension one less rules out the first possibility as well. Therefore, we get the following inductive statement towards the Campana--Peternell conjecture.

    \begin{coro} \label{coro:inductive}
		Rigidity of the canonical class for smooth projective varieties with Kodaira dimension zero in dimension $n-1$ and non-vanishing in dimension $n$ imply the Campana--Peternell conjecture in dimension $n$.
	\end{coro}

    \begin{proof}
        By \cite{Schnell-conj}*{\S9}, the Campana--Peternell conjecture on dimension $n$ reduces to non-vanishing in dimension $\leq n$ and Conjecture \ref{conj:alg-fibr-space} for $\dim X = n$. Theorem \ref{theo:main-theorem} (2) shows that Conjecture \ref{conj:alg-fibr-space} is true under the assumption that the general fiber has a rigid canonical class.
    \end{proof}

    We also give an unconditional result for fourfolds.
    \begin{coro} \label{coro:fourfolds}
        Let $X$ be a smooth projective fourfold. Let
        $$ \kappa_{\mathrm{CP}}(X) = \sup \{ \kappa(X, D) : mK_{X} - D \text{  is pseudo-effective for some } m > 0\}.$$
        Then $\kappa(X) = \kappa_{\mathrm{CP}}(X)$ unless $(\kappa(X), \kappa_{\mathrm{CP}}(X)) = (-\infty, 0)$.
    \end{coro}
    \begin{proof}
        It is clear that $\kappa_{\mathrm{CP}}(X) \geq \kappa(X)$, hence we assume that $\kappa_{\mathrm{CP}}(X) \geq 1$. Let $D$ be a divisor on $X$ such that $mK_{X} - D$ is pseudo-effective. Then following \cite{Schnell-conj}, after a birational modification, we can assume that we have an algebraic fiber space $f \colon X \to Y$ such that $D = f\sta H$ for an ample divisor $H$ on $Y$. Since the canonical bundle of the general fiber is pseudo-effective, a higher multiple is effective, since abundance for threefolds is known \cite{Kaw-abundance-threefolds}. Hence, we can apply the strategy in \cite{Schnell-conj}*{\S9} to reduce to the case when the general fiber of $f$ has Kodaira dimension zero, after suitably changing $D$. Then we apply Theorem \ref{theo:main-theorem}(2).
    \end{proof}

    \begin{rema}
        We remark that Corollary \ref{coro:fourfolds} is slightly stronger than the inductive formulation in Corollary \ref{coro:inductive} since we didn't assume the non-vanishing conjecture in dimension 4, but rather treated it as an `exceptional case'. We point out that Corollary \ref{coro:fourfolds} can already be achieved using the method by Taji that is explained above.
    \end{rema}

    Finally, we mention that the Campana--Peternell conjecture is closely related to the `birational stability conjecture' of Campana (see \cite{CP-conj}*{page 43 line 3}). For a smooth projective variety $X$ of dimension $n$, we define
    $$ \kappa_{+}(X) = \sup \{ \kappa(\det \cF) : \cF \subset \Omega_{X}^{p} \text{ for some } p > 0\}.$$
    The conjecture predicts that $\kappa_{+}(X) =\kappa(X)$ if $K_X$ is pseudo-effective. This is actually a special case of the Campana--Peternell conjecture as mentioned in \cite{CP-conj}. Indeed, using \cite{Campana-Paun}*{Theorem 1.3} or \cite{CP-conj}*{Theorem 1.8}, $L = \det \cF$ embeds to $\otimes^{m} \Omega_{X}^{1}$ for some $m > 0$ and we see that
    $$ \det \left( \otimes^{m} \Omega_{X}^{1} / L\right) \simeq \omega_{X}^{ \otimes m \cdot n^{m-1}} \otimes L^{-1}$$
    is pseudo-effective (one can assume that $L$ is saturated in $\otimes^{m} \Omega_X^1$). Therefore, the Campana--Peternell conjecture tells $\kappa(X) \geq \kappa(L)$. The other inequality $\kappa(X) \leq \kappa_{+}(X)$ is clear.

    \section{The proof}
	\subsection{Strategy of the proof} We give an outline of the proof of Theorem \ref{theo:main-theorem}. The main strategy is to use the Siu decomposition \cite{Siu-decomp} to deal with the divisorial part of closed positive (1,1)-currents. It says that for a closed positive (1,1)-current $T$ on a complex manifold $X$, $T$ can be written as
	$$ T = \sum_{W} \nu(T, W)[W] + R, $$
	where $W$ runs through all codimension 1 analytic subsets of $X$, $\nu(T, W)$ is the generic Lelong number of $T$ along $W$, and $R$ is a closed positive (1,1)-current such that $\nu(R, W) = 0$ for all codimension 1 analytic subsets of $X$. Also, the sum in the formula above is a countable sum. Lelong numbers could be understood as analytic analogues of multiplicities and we refer to \cite{Dem-anal-methods}*{Chapter 2} for the basic definition and properties of Lelong numbers. The generic Lelong number is defined as
	$$ \nu(T, W) = \inf_{x \in W} \nu(T, x)$$
	and the equality holds for a very general member $x \in W$ due to the analyticity of the super-level sets of Lelong numbers (see \cite{Siu-decomp} or \cite{Dem-anal-methods}*{14.3}).
	
	The way we use the Siu decomposition is the following. Let $L$ be a pseudo-effective divisor class on $X$, hence equipped with a singular metric $h_{L}$ with semi-positive curvature. Denote by $T_{h_{L}}$ the curvature current of $h_{L}$. Then for a prime divisor $W$, the divisor class $L - \nu(T_{h_{L}}, W)[W]$ is again pseudo-effective. Hence, we can subtract certain effective divisors from pseudo-effective classes and still remain pseudo-effective.

	Of course, we have to start with the singular metric on the line bundle $L_{0} = m_{0}K_{X} - f\sta H$, where we have a priori no control on the singularity of the metric. Under the assumption (1) in Theorem \ref{theo:main-theorem}, we use the Bergman-kernel metric on
	$$L_{1} = m_{1}K_{X/Y} + (m_{0}K_{X} - f\sta H)$$
	due to \cite{Paun-Takayama}. The point in \cite{Schnell-conj} is that the singularities of the metric improve via considering the Bergman-kernel metric on $L_{1}$. We want to stress that when the Kodaira dimension of the general fiber is zero, the description of the Bergman-kernel metric is even more explicit. However, due to the presence of the \textit{relative} pluricanonical bundle on $L_{1}$, Schnell has to assume that $K_{Y}$ is pseudo-effective in order to prove \cite{Schnell-conj}*{Theorem 12.1}. Similarly, we cannot avoid this issue as well, but we get a weaker assumption in this case by explicitly analyzing the singularity of the Bergman-kernel metric over the codimension 1 locus of $Y$. As a result, we gain an extra effective factor of $(1-\eps)B_{Y}$.

    Under the second assumption, we observe that the rigidity of the canonical class of the fibers gives a strict constraint on the behavior of a positively curved metric $h_{L_{0}}$ on $L_{0} = m_{0} K_{X} - f\sta H$, since the restrictions of $h_{L_{0}}$ to the fibers are essentially unique up to scaling.

    \subsection{Coefficients of $\Delta$ and $B_{Y}$ in the canonical bundle formula} \label{subsec:coeff-CBF}
	We describe some coefficients of $\Delta$ appearing in the canonical bundle formula (\ref{equation:CBF}) in \S\ref{subsec:CBF}. Let $E$ be a horizontal divisor on $X$. Then we have
	$$ r\ord_{E}\Delta = \ord_{E|_{F}} \sigma $$
	where $\sigma \in H^{0}(F, rK_{F})$ is the unique (up to constant) section. Here we take a sufficiently divisible $r > 0$.
	
	We describe the coefficients of the vertical divisors lying over codimension 1 points of $Y$. Consider $Q_{l} \subset Y$ and let
	$$ J_{l} := \{ j \in J : f(P_{j}) = Q_{l}\}.$$
	We consider a local non-vanishing section $\alpha$ of $f\lsta \omega_{X/Y}^{\otimes r}$ in a neighborhood $\Omega$ of a general point of $Q_{l}$. By fixing a local volume form $dt$ on $\Omega$, we obtain an $r$-pluricanonical form $\alpha(dt^{\otimes r})$ on $f^{-1}(\Omega)$. For $j \in J_{l}$, the order of vanishing of $\alpha(dt^{\otimes r})$ along $P_{j}$ equals to $r \cdot \ord_{P_{j}}\Delta$. One can also see that for all $l \in L$, we have $\ord_{P_{j}}\Delta \geq 0$ for $j \in J_{l}$, and there exists at least one $j \in J_{l}$ such that $\ord_{P_{j}}\Delta = 0$. We put $a_{j} = \ord_{P_{j}}\Delta$.
	
	We now explain how the coefficients of $B_{Y}$ are determined. Let $B_{Y} = \sum b_{l} Q_{l}$. Then,
	$$ 1 - b_{l} = \sup \{ b \in \QQ : (X, -\Delta + b\cdot f\sta Q_{l}) \text{ is lc over } \eta_{Q_{l}}\}$$
	where $\eta_{Q_{l}}$ is the generic point of $Q_{l}$. These numbers can be described more explicitly as follows. Put $\delta_{j} = \ord_{P_{j}} f\sta Q_{l}$ for $j \in J_{l}$. Then we can easily see that
	$$  b_{l} = \max_{j \in J_{l}} \left( 1- \frac{a_{j}+1}{\delta_{j}} \right)$$
	and that $b_{l}$ lies in the interval $[0, 1[$. We also refer to \cite{Kollar-KodCBF} for a nice exposition on the canonical bundle formula.

    \subsection{The proof for the second case} 
	We give a proof of Theorem \ref{theo:main-theorem} under the second assumption. We recall that $L_{0} = m_{0}K_{X} - f\sta H$ is pseudo-effective, hence carries a semi-positively curved metric $h_{0}$. Throughout the paper, we will always take $m_{0}$ to be divisible enough. We first show that we can subtract the horizontal part of $\Delta$.

    \begin{lemm}
        Assume that the canonical class of the general fiber $K_{F}$ is rigid, as in the second assumption of the setup of Theorem \ref{theo:main-theorem}. Then, $m_{0}(K_{X} - \Delta^{h}) - f\sta H$ is pseudo-effective.
    \end{lemm}

    \begin{proof}
        Let $Y_{\mathrm{rig}}$ be the set of points $y \in Y_{0}$ such that the canonical class of the fiber $K_{X_{y}}$ is rigid. By assumption, $Y_{\mathrm{rig}}$ contains a Zariski open subset of $Y_{0}$. Note that $L_{0}|_{X_{y}} = m_{0}K_{X_{y}}$. Hence, for each $y \in Y_{\mathrm{rig}}$, the restricted metric $h_{0}|_{X_{y}}$ is either constantly $+\infty$, or its curvature current is $m_{0}[\Delta^{h}|_{X_{y}}]$. We denote by $T_{h_{0}}$ the curvature current of $h_{0}$. Consider an open ball $\Omega \subset Y_{\mathrm{rig}}$ and consider a section $\sigma$ on $H^{0}(f^{-1}(\Omega), m_{0}K_{X} - f\sta H)$ which vanishes along $m_{0}\Delta^{h}$. We see that $|\sigma|_{h_{0}} \colon f^{-1}(\Omega) \to \RR \cup \{ +\infty\}$ is constant along the fibers. 
        
        Let $E$ be one of the components $\Delta^{h}$ and let $\alpha_{E} = \ord_{E} \Delta^{h}$. Then $\sigma$ vanishes along $E$ with order $m_{0}\alpha_{E}$. Consider an open subset $\Omega' \subset f^{-1}(\Omega)$ intersecting $E$. We may assume that $\Omega'$ does not intersect other components of $\Delta^{h}$. Suppose that $E$ is given by the equation $u = 0$. Then $\tau \coloneqq \sigma / u^{m_{0}\alpha_{E}}$ is a local non-vanishing section of $m_{0}K_{X} - f\sta H$ on $\Omega'$. By the rigidity assumption, we have
        $$ |\tau|_{h_{0}}^{2} = e^{-\varphi(y)} |u|^{-2m_{0} \alpha_{E}}.$$
        Note that $\varphi$ only depends on $y \in \Omega$ and it is psh on $\Omega$. In particular, it is bounded above. This shows that $\nu(T_{h_{0}}, E) \geq m_{0} \alpha_{E}$ since $E \cap \Omega'$ is Zariski dense in $E$. Using the Siu decomposition, we get the pseudo-effectivity of $m_{0}(K_{X} - \Delta^{h}) - f\sta H$.
    \end{proof}

    Now we give the proof of the theorem.

    \begin{proof}[Proof of Theorem \ref{theo:main-theorem} (2)]
    Let $L = K_{Y} + B_{Y} + M_{Y}$. Note that over $Y_{0}$, we have
	$$ m_{0}(K_{X} - \Delta^{h}) -f\sta H \sim f\sta (m_{0}L - H).$$
	For a local non-vanishing section $\sigma \in H^{0}(\Omega, m_{0}L - H)$ for $\Omega \subset Y_{0}$, we see that
	$$ |\sigma|_{\widetilde{h}_{0}, x}^{2} : X_{y} \to \RR_{\geq 0} \cup \{+\infty\} $$
	has to be a constant function since the only global psh functions on a compact space are constants. We define $|\sigma|_{h, y} := |\sigma|_{\widetilde{h}_{0}, x}$ for $y = f(x)$. This gives a positively curved singular metric on $m_{0}L - H$ over $Y_{0}$. Therefore, it remains to show that the local potentials are bounded above outside of a codimension 2 subset of $Y$. Then this metric automatically extends across $Y \setminus Y_{0}$ and therefore, $m_{0}L - H$ is pseudo-effective.
	
	We recall from \S\ref{subsec:CBF} for the set-up of the canonical bundle formula and \S\ref{subsec:coeff-CBF} for the description of the coefficients appearing therein. Consider an open neighborhood $\Omega$ of a point of $Q_{l}$ that avoids $f(\supp \Delta^{-})$ and other components of $Q$. Let $\sigma$ be a local non-vanishing section of $m_{0}L - H$ on $\Omega$, which we can also view as a section on $m_{0}(K_{X} - \Delta) - f\sta H$ on $f^{-1}(\Omega)$, vanishing along the divisor $\sum_{j \in J_{l}} m_{0} a_{j} P_{j}$. Note that there exists some $j \in J_{l}$ such that $a_{j} = 0$. Then we have an open set $\Omega'$ of $f^{-1}(\Omega)$ intersecting $P_{j}$ that surjects to $\Omega$. Hence, we can assume that $\sigma$ does not vanish on $\Omega'$. Since $|\sigma|_{h_{0}, x}^{2}$ is bounded below by a positive constant on $\Omega'$, $|\sigma|_{\widetilde{h}_{0}, y}^{2}$ is also bounded below by a positive constant on $\Omega$. This shows that the local potentials for the metric $\widetilde{h}_{0}$ are bounded above. This concludes the proof of the second case. 
    \end{proof}

    The rest of the paper is devoted to the proof of the first case of Theorem \ref{theo:main-theorem}.

	\subsection{Bergman-kernel metric} We briefly describe the construction of the relative Bergman-kernel metric, constructed in \cite{Paun-Takayama}, adapted to our setup. The upshot is that for a given metric $h_{0}$ on $L_{0} = m_{0}K_{X} - f\sta H$, we can induce a positively curved hermitian metric on $L_{1}= m_{1} K_{X/Y} + L_{0}$ with `better properties', when $m_{1}$ is big enough. We point out that the situation is substantially simpler than \cite{Paun-Takayama} due to the fact that the Kodaira dimension of the general fiber is zero. We assume that $m_{0}$ and $m_{1}$ are sufficiently divisible and $m_{1} \gg m_{0}$. We let $f^{-1}(Y_{0}) = X_{0}$ and describe the singular hermitian metric of $h_{1}$ on $L_{1} = m_{1} K_{X/Y} + L_{0}$ on $X_{0}$. For $y \in Y_{0}$ and $\sigma_{y} \in H^{0}(X_{y}, m_{1}K_{X_{y}} + L_{0,y})$, we define
	$$ l(\sigma_{y}) := \left( \int_{X_{y}} |\sigma_{y}|^{2/m_{1}}_{h_{0}}\right)^{m_{1}/2}.$$
	For the right-hand side, we can view $\sigma_{y}$ as an $m_{1}$-pluricanonical form with values in $(L_{0, y}, h_{0, y})$ on $X_{y}$, and therefore, $|\sigma_{y}|^{2/m_{1}}_{h_{0}}$ is a well-defined volume form on $X_{y}$. We simply define $l(\sigma_{y}) = \infty$ when $h_{0, y} \equiv +\infty$. For a fixed point $x \in X_{0}$ such that $f(x) = y$ and a local section $s$ of $L_{1}$ at $x$, we have
	$$ |s|_{h_{1}, x} = l(\sigma_{y}),$$
	where $\sigma_{y}(x) = s(x)$. This is due to the fact that the dimension of $H^{0}(X_{y} , m_{1}K_{X_{y}} + L_{0,y})$ is 1, so the supremum in the formula goes away. We follow the convention that if $\sigma_{y}$ vanishes at $x$, then $|s|_{h_{1},x} = +\infty$. The main result in \cite{Paun-Takayama} is that the local potentials of this metric are upper semi-continuous and plurisubharmonic, and bounded above towards $X \setminus X_{0}$. This allows us to extend the metric on $L_{1}|_{X_{0}}$ uniquely to a semi-positively curved singular metric on $L_{1}$. We note that since $m_{1} \gg m_{0}$, the set $\{ y \in Y_{0} : l(\sigma_{y}) < +\infty\}$ is non-empty, and in fact, of full measure. We denote by $T_{1}$ the curvature current associated to the metric $h_{1}$ on $L_{1}$.
	
	\begin{rema}
		Following the notation in \cite{Paun-Takayama}, one technical issue in the construction of the Bergman-kernel metric is that they have to define the set
		$$ Y_{m_{1}, \mathrm{ext}} = \{ y \in Y_{0} : h^{0}(X_{y}, m_{1}K_{X_{y}} + L_{0, y}) = \rank f\lsta (m_{1}K_{X/Y} + L)\}$$
		and construct the metric on $X_{m_{1}, \mathrm{ext}} = f^{-1}(Y_{m_{1}, \mathrm{ext}})$, ignoring the behavior of the metric on $X_{0} \setminus X_{m_{1}, \mathrm{ext}}$. In our situation, we have $Y_{0}  = Y_{m_{1}, \mathrm{ext}}$ due to the invariance of plurigenera \cite{Siu-two-tower, Paun-one-tower}.
	\end{rema}
	
	\subsection{Lelong numbers along the horizontal divisors}
	We first estimate the Lelong numbers on the horizontal divisors. Let $E$ be a component of $\Delta^{h}$ and let $\alpha = \ord_{E} \Delta^{h}$. Consider an open set $\Omega \subset Y_{0}$ and fix a local non-vanishing section $\sigma \in f\lsta (m_{1}K_{X/Y} + L_{0})$. We may assume that $1\leq l(\sigma_{y}) < +\infty$ for $y \in \Omega$ by shrinking $\Omega$ and rescaling $\sigma$. Let $s$ be a local non-vanishing section of $m_{1}K_{X/Y} + L_{0}$ on a neighborhood of $E \cap f^{-1}(\Omega)$ and suppose that $E$ is given by the equation $u = 0$. Then we can assume that $\sigma =  u^{\alpha(m_{1} + m_{0})} s$, where we consider $\sigma$ as an element in $H^{0}(f^{-1}(\Omega), m_{1}K_{X/Y} + L_{0})$. We denote by $|s|_{h_{1}}^{2} = e^{-\varphi}$ and we get
	$$ e^{-\varphi} =  |u|^{-2\alpha(m_0 +m_1)}l(\sigma_{f(x)})^{2} \geq |u|^{-2\alpha(m_{0}+m_{1})}.$$
	This implies that $\nu(\varphi, x) \geq \alpha(m_{1}+m_{0})$ for $x \in E \cap f^{-1}(\Omega)$. Therefore, we have $\nu(T_{1}, E) \geq \alpha(m_{1}+m_{0})$ since $E \cap f^{-1}(\Omega)$ is Zariski dense in $E$.
	
	\subsection{Lelong numbers along the vertical divisors} \label{subsec:Lelong-vertical} Here, we fix $Q_{l}$ and $P_{j}$ such that $j \in J_{l}$. We recall the notations for the setup from \S\ref{subsec:coeff-CBF}. Fix a general point $x_{0} \in P_{j}$ so that $x_{0}$ is not contained in other divisors in the support of $\Delta$. We fix local coordinates $t_{1},\ldots, t_{p}$ on $Y$ centered at $f(x_{0})$ and $z_{1},\ldots, z_{q}$ on $X$ centered at $x_{0}$ such that 
	\begin{itemize}
		\item $P_{j}$ is defined by the equation $z_{1} = 0$,
		\item $Q_{l}$ is defined by the equation $t_{1} = 0$, and
		\item The local expression of $f:X\to Y$ is $t_{1} = z_{1}^{\delta_{j}}$, and $t_{j} = z_{j}$ for $j = 2, \ldots, p$.
	\end{itemize}
	Let $\sigma$ be a local non-vanishing section of $f\lsta (m_{1}K_{X/Y} + L_{0})$ on a neighborhood of $f(x_{0})$, and let $s$ be a non-vanishing section of $m_{1}K_{X/Y} + L_{0}$ on a neighborhood of $x_{0}$. Then by the discussion in \S\ref{subsec:coeff-CBF}, we may assume that
	$$ \sigma =  z_{1}^{(m_{1}+m_{0})a_{j}} s$$
	and therefore
	$$ |s|_{h_{1}, x}^{2} = |z_{1}|^{-2(m_{1}+m_{0}) a_{j}} l(\sigma_{f(x)})^{2}. $$
	Using direct calculation, we will show separately in \S\ref{section:vol-asymptotics} that for fixed $\eps > 0$, we have
	$$ l(\sigma_{f(x)})^{2} \gtrsim |t_{1}|^{-2m_{1}(1-\eps)b_{l}},$$
	if $m_{1}$ is sufficiently large. The proof of the volume asymptotics will be omitted at this point in order to simplify the arguments and stress the key ideas. Since $t_{1} = z_{1}^{\delta_{j}}$, we have that
	$$ \nu(T_{1}, P_{j}) \geq (m_{1}+m_{0})a_{j} + (1-\eps)b_{l}\delta_{j} m_{1}.$$
	
	\subsection{Conclusion}
	From the computation of the generic Lelong numbers of the metric on $L_{1}$, we see that
	$$L_{1}' = m_{1}K_{X/Y} + (m_{0}K_{X} - f\sta H) - (m_{1} + m_{0})\left(\Delta^{h} + \sum_{l \in L}\sum_{j \in J_{l}} a_{j}P_{j}\right) - m_{1}\sum_{l \in L}\sum_{j \in J_{l}} (1-\eps)b_{l}\delta_{j} P_{j} $$
	is pseudo-effective. Over a complement of a codimension 2 subset of $Y$, we see that this is equivalent to
	$$ f\sta \Big((m_{1}+m_{0})(K_{Y} +B_{Y} + M_{Y}) - m_{1}(K_{Y} + (1-\eps)B_{Y}) - H\Big).$$
	Hence, the singular hermitian metric on $L_{1}'$ descends to a semi-positively curved singular hermitian metric on
	$$(m_{1} +m_{0})(K_{Y} + B_{Y} + M_{Y}) - m_{1}(K_{Y} + (1-\eps)B_{Y}) - H,$$
	over a complement of a codimension 2 subset of $Y$. Therefore, the metric extends to the whole space $Y$. Since we assumed that $K_{Y} + (1-\eps)B_{Y}$ is pseudo-effective, $m_{1}(K_{Y} +(1-\eps)B_{Y}) + H$ is big, and therefore, $K_{Y} + B_{Y} + M_{Y}$ is big. This concludes the first case of Theorem \ref{theo:main-theorem}. \hfill{$\square$}
	
	\begin{rema}
		The argument unconditionally implies that $M_{Y} + \delta(K_{Y} + B_{Y}) + \eps B_{Y}$ is big for some $0 < \delta \ll \eps \ll 1$.
	\end{rema}
	
	\begin{ques}
		Even though $K_{Y} + (1-\eps)B_{Y}$ is not pseudo-effective in general, we predict that $K_{Y} + Q$ is big. We already know that $M_{Y} + \delta (K_{Y} + B_{Y}) + \eps B_{Y}$ is big for some $0 < \delta \ll \eps \ll 1$. Indeed, if one manages to show that $M_{Y}$ is big, then $K_{Y} + Q$ is big. 
        
        Let's assume that $M_{Y}$ is big. Note from \cite{Kollar-KodCBF} that we have a variation of Hodge structures $\VV$ on $Y_{0}$ and after a finite cover $\pi \colon \widetilde{Y} \to Y$, the pull-back of $M_{Y}$ coincides with the canonical extension of the lowest piece of $\pi\sta \VV$ and we can also assume that $\pi\sta \VV$ has unipotent local monodromies. Then one can calculate the intersection number $(\pi\sta M_{Y})^{\dim Y}$ by integrating the curvature induced by the polarization of $\VV$, and this being strictly positive implies that the period map is immersive at one point (see \cite{BBT-ominimal-GAGA}*{Lemma 6.15 and 6.17} for instance). This shows that the period map on $Y_{0}$ is immersive at least at one point. Then the result of \cite{Brunebarbe-Cadorel} says that $K_{Y} + Q$ is big. 

        In this vein, we would like to ask the question if one can unconditionally show that $K_{Y} + Q$ is big in the set-up of Theorem \ref{theo:main-theorem}. A posteriori, this should be the case since we have $\kappa(X) = \dim Y$ and this implies the bigness of $K_{Y} + Q$ by \cite{park2022logarithmic}*{Theorem 1.9}.
	\end{ques}

	\section{Volume asymptotics} \label{section:vol-asymptotics}
	Here, we give the proof of the volume asymptotics in \S\ref{subsec:Lelong-vertical}.
    \begin{lemm}
        As in the setup of \S\ref{subsec:Lelong-vertical}, we have the following volume asymptotics:
        $$ l(\sigma_{f(x)})^{2} \gtrsim |t_{1}|^{-2m_{1}(1-\eps)b_{l}},$$
        if $m_{1}$ is sufficiently big.
    \end{lemm}

    \begin{proof}
        First, we point out that we can ignore the case when $b_{l} =0$ since $l(\sigma_{f(x)})$ is bounded below by a positive constant following the proof of \cite{Paun-Takayama}*{Theorem 4.2.7}. Therefore, we assume that $b_{l}> 0$.
	
	Let $\Omega$ be a neighborhood on a general point of $Q_{l}$ and write the coordinate functions as $(t_{1},\ldots, t_{p})$ and identify $\Omega$ with the polydisk of radius 1. Assume that $Q_{l}$ is defined by the equation $t_{1} = 0$. Let
	$$ f^{-1}(\Omega) \supset \bigcup_{j \in J_{l}} \Omega\uind{j},$$
	where $\Omega\uind{j}$ are disjoint open sets where $\Omega\uind{j}$ only intersects with $P_{j}$. Fix coordinate charts $z_{1}\uind{j},\ldots, z_{q}\uind{j}$ of $\Omega\uind{j}$ and assume that $P_{j}$ is defined by the equation $z_{1}\uind{j}$ and the local expression for $f$ is
	$$ f(z_{1}\uind{j},\ldots, z_{q}\uind{j}) = ((z_{1}\uind{j})^{\delta_{j}}, z_{2}\uind{j}, \ldots, z_{p}\uind{j}). $$
	Identify $\Omega\uind{j}$ with the polydisk of radius 1 as well. Let $\sigma$ be a local non-vanishing section of $f\lsta (m_{1}K_{X/Y} + L_{0})$ on $\Omega$ and $s_{j}$ be local non-vanishing sections of $m_{1}K_{X/Y} + L_{0}$ on $\Omega\uind{j}$. We can assume that
	$$ \sigma = (z_{1}\uind{j})^{(m_{1} + m_{0})a_{j}} s_{j}$$
	and
	$$ s_{j}(dt^{\otimes m_{1}}) = (dz_{1}\uind{j}\wedge \ldots \wedge dz_{q}\uind{j})^{\otimes m_{1}} \otimes e_{j},$$
	where $e_{j}$ is a non-vanishing section of $L_{0}$ on $\Omega\uind{j}$. Here, we identify $m_{1}K_{X/Y} + L_{0}$ with the sheaf
	$$ \SheafHom_{\cO_{X}} (f\sta (m_{1}K_{Y}) , m_{1}K_{X} + L_{0}).$$
	We let $|e_{j}|_{h_{0}}^{2} = e^{-\varphi_{j}}$, where $\varphi_{j} : \Omega\uind{j} \to \RR \cup\{-\infty\}$ are the local weight functions for the metric on $L_{0}$. By shrinking the domains if necessary, one can assume that $\varphi_{j} \leq C_{0}$ for some constant. For $y = (t_{1},\ldots, t_{p}) \in \Omega$, we have
	$$ l(\sigma_{y})^{2/m_{1}} = \int_{X_{y}} |\sigma_{y}|_{h_{0}}^{2/m_{1}} \geq \sum_{j \in J_{l}} \int_{\Omega_{y}\uind{j}} |\sigma_{y}|_{h_{0}}^{2/m_{1}},$$
	where $\Omega_{y}\uind{j} = X_{y} \cap f^{-1}(\Omega)$. We note that
	$$ \frac{dz_{1}\uind{j}}{z_{1}\uind{j}} \wedge dz_{2}\uind{j} \wedge \ldots \wedge dz_{p}\uind{j} = \frac{1}{\delta_{j}} f\sta \left(\frac{dt_{1}}{t_{1}} \wedge dt_{2} \wedge \ldots \wedge dt_{p}\right).$$
	By definition, $\sigma_{y}$ is the section on $m_{1}K_{X_{y}}+L_{0, y}$ satisfying
	$$ \sigma(dt^{\otimes m_{1}})|_{X_{y}} = \sigma_{y} \otimes f\sta (dt^{\otimes m_{1}}).$$
	Therefore, the local expression for $\sigma_{y}$ is
	$$ \sigma_{y} = \frac{1}{\delta_{j}^{m_{1}}} (z_{1}\uind{j})^{(m_{1}+m_{0})a_{j} + (1-\delta_{j})m_{1}} \cdot (dz_{p+1}\uind{j} \wedge \ldots \wedge dz_{q}\uind{j})^{\otimes m_{1}} \otimes s_{j}.$$
	This shows that
	$$ \int_{\Omega_{j}\uind{j}} |\sigma_{y}|_{h_{0}}^{2/m_{1}} \geq C_{1} \cdot |t_{1}|^{2 \cdot \frac{(m_{1} + m_{0})a_{j} + (1-\delta_{j})m_{1}}{\delta_{j}m_{1}}},$$
	where $C_{1}>0$ is some constant, using $t_{1} = (z_{1}\uind{j})^{\delta_{j}}$.
	Therefore, we can see that
	$$ \int_{X_{y}} |\sigma_{y}|_{h_{0}}^{2/m_{1}} \gtrsim |t_{1}|^{-2b_{l}'},$$
	where
	$$ -b_{l}' = \min_{j \in J_{l}} \frac{(m_{1} + m_{0}) a_{j} + (1-\delta_{j}) m_{1}}{m_{1}\delta_{j}} \geq \min_{j \in J_{l}} \frac{m_{1} a_{j} + (1-\delta_{j}) m_{1}}{m_{1}\delta_{j}} = -b_{l}.$$
	Moreover, we see that
	$$ |b_{l} - b_{l}'| \leq \frac{m_{0}}{m_{1}} \max_{j \in J_{l}} \frac{a_{j}}{\delta_{j}}$$
	and we can make the right-hand side to be less than $\eps b_{l}$, by taking $m_{1}$ sufficiently big, since $b_{l} > 0$. This shows that
	$$ l(\sigma_{y})^{2} \gtrsim |t_{1}|^{-2m_{1}(1-\eps)b_{l}}.$$
    \end{proof}
	
	\begin{rema}
		If one uses the locus on the intersection of $P_{j}$'s, it is possible to get an extra log term for the volume asymptotics, following the computation of \cite{Boucksom-Jonsson} or \cite{Takayama-singularitiesNS}. However, we do not write it down since we don't need this.
	\end{rema}
	
	\textbf{Acknowledgements.} The author would first like to thank the organizers of the conference `Transcendental aspects of algebraic geometry' in Cetraro 2024 where the author was able to attend a lovely talk by Christian Schnell on this topic. The author would also like to thank his advisor Mircea Musta\c{t}\u{a} for generous support, Mattias Jonsson, Vladimir Lazi\'{c}, Niklas Müller, Sung Gi Park, Mihnea Popa, Mihai P\u{a}un for useful conversations during the preparation of this article, and also the referees for helpful suggestions.
	
	\bibliographystyle{alpha}
	\bibliography{Reference}

\end{document}